\documentclass[12pt,a4paper]{article}
\def\R{\mathbf{R}}
\usepackage{amsmath}
\usepackage{amsfonts}
\usepackage{amssymb}
\usepackage{xspace}
\usepackage{graphics}

\renewcommand\iff{if and only if\xspace}

\def\l{\left}  \def\r{\right}
\def\a{\alpha} 
\def\const{\mathrm{const}}
\def\Ima{\mathrm{Im}\ }
\def\Rea{\mathrm{Re}\ }
\renewcommand\t{\tau}

\newcommand\LB[1]{\label{#1}} 
\newcommand\BE[2]{\begin{#1} #2 \end{#1}}

\newcommand\EQ[2]{\BE{equation}{\LB{#1} #2}}

\newcommand{\eq}{equation\xspace}

\newcommand\ST{such that }
\newcommand{\lint}{\int\limits}

\newcommand\f{\varphi}

\newcommand\gm{\gamma}

\newcommand \lm{\lambda}
\newcommand \Lm{\Lambda}

\newcommand\dl{\delta}

\newcommand{\bbC}{\blackboard{C}}

\newcommand{\bbR}{\blackboard{R}}
\newcommand{\blackboard}[1]{\mathbb#1}
\newcommand{\iy}{\infty}

\title{Spectra of the Gurtin-Pipkin type equations
}
\author{A. Eremenko\thanks{Supported by NSF
and by the Humboldt Foundation.}$\;$
and
S. Ivanov\thanks{
Supported in part by the Russia Foundation for Basic
Research, grant 08-01-00595a.
}}

\begin{document}

\maketitle
\begin{abstract}
We study the spectra of certain
integro-differential equations arising in applications.
Under some conditions on the kernel of the integral operator,
we describe the non-real part of the spectrum.

MSC Primary: 45K05, 35P20. Keywords: Gurtin--Pipkin equation, spectrum,
Denjoy--Wolff theorem.

\end{abstract}

\section{\LB{intro_etc} Introduction}


The following integro-differential equations
arise in several areas of physics and applied mathematics,
namely in
heat transfer with finite propagation speed \cite{GuPip},
systems with thermal memory \cite{EMV},
viscoelasticity problems \cite{CMD}
and acoustic waves in composite media \cite{Sham}.

(i) Gurtin--Pipkin equations of first order in time (GP1)
\begin{equation}\LB{1GP1}
u_{t}(x,t)=\int_0^t k(t-s)
u_{xx}(x,s)\,ds, \ x\in(0,\pi), \ t>0,
\end{equation}

(ii) and of second order in time (GP2)

\begin{equation}\LB{2GP2}
u_{tt}(x,t)=a u_{xx}-\int_0^t k(t-s) u_{xx}(x,s)\,ds,\ a>0,\ x\in(0,\pi),  \ t>0.
\end{equation}

(iii)
Kelvin-Voigt equation (KV)
\EQ{KV}{
u_{tt}(x,t)=u_{xx}+\epsilon u_{txx}-\int_0^t k(t-s) u_{xx}(x,s)\,ds, \ x\in(0,\pi), \ t>0.
}

Here
$$
k(t)=\int_0^\infty e^{-t\t}d\mu(\t)
$$
is the Laplace transform of a positive measure $d\mu$.
We identify this measure with its distribution
function $\mu$, so $\mu$ is increasing,
continuous from the right, and the integral
is interpreted as a Stieltjes integral \cite{W}.
We always assume that
$k$ is defined and integrable on $(0,\infty)$, that is
\begin{equation}\label{1}
\int_0^\infty\frac{d\mu(t)}{t}<\infty,
\end{equation}
and that $\mu$ is supported on $(d_0,\infty)$ with
some $d_0>0$.

\BE{remark}{If equation \eqref{1GP1} can be differentiated
with respect to $t$,
then we obtain a special
case of \eqref{2GP2}:
$$
u_{tt}=k(0)u_{xx}-\int_0^t \tilde k(t-s)u_{xx}ds.
$$
with
\EQ{dk1}{
\tilde  k=-\frac d{dt}k.
}
}

\BE{remark}{
In the case $k(t)=\const=\a^2$, \eq \eqref{1GP1} is
in fact an integrated wave \eq. Indeed,
differentiation of \eqref{1GP1}  gives
\EQ{Wave}{
u_{tt}=\a^2 u_{xx}.
}
If $k(t)=\a^2e^{-bt}$,
then differentiation gives a damped wave \eq
\EQ{wave1}{
u_{tt}=\a^2u_{xx}-bu_t.
}
}

Let the initial conditions be
$u(\cdot,0)=\xi$ for \eqref{1GP1}, and $u(\cdot,0)=\xi$,  $ u_t(\cdot,0)=\eta$
for \eqref{2GP2} and \eqref{KV}.

First we apply Fourier's method:
we set $\f_n=\sqrt{\frac2\pi}\sin nx$ and expand the solution and the initial
data in a series in $\f_n$
$$
u(x,t)=\sum_1^\iy u_n(t)\f_n(x),\ \xi(x)=\sum_1^\iy \xi_n \f_n(x),\
 \eta(x)=\sum_1^\iy \eta_n \f_n(x).
$$

The components $u_n$ satisfy ordinary
integro-differential equations
\vspace{.1in}

(i) GP1
\begin{equation}\label{thn1}
\dot u_n(t)=-n^2\int_0^t k(t-s) u_n (s)d s, \  t>0.
\end{equation}

(ii) GP2
\begin{equation}\label{thn2}
\ddot u_n(t)=-a n^2u_n(t)+n^2\int_0^t k(t-s) u_n (s)d s, \  t>0.
\end{equation}

(iii) KV
\begin{equation}\label{thnKV}
\ddot u_n(t)=- n^2u_n(t)-\epsilon n^2\dot u_n+n^2\int_0^t k(t-s) u_n (s)d s, \  t>0.
\end{equation}

We will denote the Laplace images by the capital letters.
Applying the Laplace Transform to \eqref{thn1}, \eqref{thn2},
and \eqref{thnKV}, and using the initial conditions, we find
\vspace{.1in}

(i) GP1
\EQ{Thn}{
zU_n(z)-\xi_n=-n^2 K(z) U_n(z)
}

or
\EQ{GP1Thn}{
U_n(z)=\frac{\xi_n}{z+n^2K(z)}.
}

(ii) GP2
$$
z^2U_n(z)-z\xi_n-\eta_n=-a n^2U_n(z)+n^2 K(z) U_n(z)
$$

or
\EQ{GP2Un}{
U_n(z)=\frac{z\xi_n+\eta_n}{z^2+an^2-n^2K(z)}.
}

(iii) KV
$$
z^2U_n(z)-z\xi_n-\eta_n=- n^2U_n(z)-\epsilon n^2\l[zU_n-\xi_n \r] +n^2 K(z) U_n(z),
$$

or
\EQ{KVUn}{
U_n(z)= \frac{\xi_n+z\eta_n-\epsilon n^2\xi_n}{z^2+\epsilon zn^2+n^2 -n^2K(z)}.
}

Denote the denominators in \eqref{GP1Thn}, \eqref{GP2Un}, and \eqref{KVUn} by $F_n(z)$, $G_n(z)$,
and $H_n(z)$ respectively:
$$
F_n(z)= z+n^2K(z), \quad G_n(z)=z^2+an^2-n^2K(z),$$
and

$$  H_n(z)=z^2+\epsilon zn^2+n^2 -n^2K(z).
$$
Let $F_n^0$, $G_n^0$, and $H_n^0$ be the sets of zeros of
$F_n(z)$, $G_n(z)$,
and $H_n(z)$ respectively. Set
$$
\Lm_{GP1}=\bigcup_1^\iy F_n^0,\ \Lm_{GP2}=\bigcup_1^\iy G_n^0,\ \Lm_{KV}=\bigcup_1^\iy H_n^0.
$$

\BE{definition}{The sets $\Lm_{GP1}$, $\Lm_{GP2}$,
and $\Lm_{KV}$
are the spectra of equations \eqref{1GP1}, \eqref{2GP2}
and \eqref{KV} respectively.
}
Study of the spectra is important for applications.
See \cite{Sham} and the references therein.

\BE{remark}{
Suppose that in the integro-differential equations
(\ref{1GP1}), (\ref{2GP2}) and (\ref{KV})
we  replace the zero lower limit in the integrals by
$-\iy$. Then $\lm$ is a point
of the spectrum of  $F_n$ or  $G_n$ or  $H_n$
\iff there is a solution of the form
$$
u_\lm(x,t) = e^{\lm t}\f_n(x).
$$
For such systems a semigroup approach to the
equations is possible,
(V.V. Vlasov, private communication).
}

\BE{remark}{If \eqref{dk1} holds, then
$\tilde K(z)=k(0)-zK(z)$.
}

We do not study here the regularity of
solutions and consider the solutions as
sequences $\{u_n(t)\}$.
Regularity of GP1  is studied in \cite{P05} and of GP2
in \cite{VW}.
In \cite{P05}, under the assumption that
$k(t)$ is twice continuously differentiable it was shown,
in particular,
that the solution $u(x,t)$ of \eqref{1GP1} is a continuous
$L^2(0,T)$-valued function.
In \cite{VW}, the conditions on the kernel $k(t)$
and the initial data are found \ST there
exists a strong solution: for a $\gm>0$
$$
\lint_0^\iy e^{-\gm t} \l[
\|u(\cdot,t)\|^2_{L^2(0,\pi)} +
\|u_{xx}(\cdot,t)\|^2_{L^2(0,\pi)}+
\|u_{tt}(\cdot,t)\|^2_{L^2(0,\pi)}
\r]dt<\iy.
$$

We single out  the important special case:
the case of discrete measure $\mu$ with atoms at
$b_k>0$ of mass $a_k>0$,
\begin{equation}
\label{case1}
K(z)=\sum_{k=1}^\infty\frac{a_k}{z+b_k},
\quad 0<b_1<\ldots\to+\infty,\quad a_k>0,
\end{equation}
Discrete measures arise in applications \cite{Sham},
where parameters $a_k$, $b_k$
are connected with auxiliary boundary value problems
arising under averaging.

\section{\LB{mainR}  {Main Results}}
In the general case the non-real
part of the spectrum is described as follows.

\BE{theorem}{\LB{complex}

(i)
For every $n$, each
set $F_n^0$, $G_n^0$, or $H_n^0$
contains at most
one point in the upper half-plane,
and this point, if exists, belongs to
the second quadrant.

(ii) For $n$ large enough, the set $G_n^0$
contains a point $z_n$  such that
\EQ{zn11}{
z_n=i\sqrt a \,\,n+o(n).
}

(iii)
If
\begin{equation}
\label{2}
A=\int_0^\infty d\mu(t)<\infty,
\end{equation}
then for $n$ large enough, the set $F_n^0$
contains a point $z_n$  such that
\EQ{zn11a}{
z_n=i\sqrt A \,\,n+o(n).
}
}

Under the additional assumptions on $K$,
we can find more precise asymptotics of the non-real
part of the spectrum.

\BE{theorem}{\LB{complex1}
Suppose that
\begin{equation}
\label{mu}
\mu(t)=bt^{\rho}+O(t^{\alpha}),\quad t\to\infty,
\end{equation}
where $0<\alpha<\rho<1.$ Then:
\newline
(i) for $n$ large enough, the set $F_n^0$
contains a zero $z_n$ of $F_n(z)$ such that
\EQ{zn12}{
z_n=\l(\frac {b\pi \rho}{\sin\pi\rho}\r)^{1/(2-\rho)}
e^{i\pi/(2-\rho)} n^{2/(2-\rho)}
\l(1+O(n^{2(\a-\rho)/(2-\rho)})\r),
}
(ii)
for $n$ large enough, the set $G_n^0$
contains a zero $z_n$ of $G_n(z)$ such that
\EQ{zn}{
z_n=i\sqrt{a}n+\frac{b\pi\rho}{2\sin\pi\rho}
a^{\rho/2-1}e^{i\pi(\rho/2-1)}n^\rho(1+o(1)).
}
}

In the case of a discrete measure with finite number of atoms
the next theorem was conjectured by
V. V. Vlasov and N. Rautiyan (private communication).

\BE{theorem}{\LB{complexKV} Let $\mu$ be a measure
with compact support, that is
\EQ{finite}{
k(t)=\int_{d_0}^d e^{-t\tau}d\mu(\tau),\quad 0<d_0<d<\infty.}
Then the set $\Lm_{KV}$
contains finite number of non-real points.}

In the case \eqref{case1} it is not hard to study the real spectrum of the
systems, see \cite{IS}\footnote{The assertions of Theorem~1
in \cite{IS} are correct only for
$n$ large enough.} and the discussion below.
The simplest versions of theorems \ref{complex}(iii)
and \ref{complex1}
are contained in \cite{Sham}, \cite{I}, \cite{IS}.
In \cite{IP}, a study of the spectrum shows
the lack of controllability of the system.

\section{\LB{main} {Proof of the Main Results and discussion}}


Our main tools are the Schwarz Lemma and the
Denjoy--Wolff Theorem (see, for example, \cite{Shapiro}).
\vspace{.1in}

\noindent
{\bf Schwarz's Lemma.} {\em Let $f$ be an
analytic function which
maps the upper half-plane $\bbC_+$ into itself. Then the
equation $f(z)=z$
has at most one solution $w$, and if such solution exists
then $|f'(w)|<1$,
unless $f$ is an elliptic fractional-linear
transformation.}
\vspace{.1in}

\noindent
{\bf Denjoy--Wolff Theorem}. {\em Let $f$ be an analytic
function which maps $\bbC_+$ into
itself,
and suppose that $f$ is not an elliptic
fractional-linear transformation. Then there exists
a unique point $w\in{\bbC_+}\cup\{\infty\}$
such that the iterates
$f^{*n}$ converge to $w$ uniformly on compact subsets
of $\bbC_+$, the angular limit $f(w)=\lim_{z\to w} f(z)$
exists and satisfies $w=f(w)$.
Moreover, the angular derivative $f'(w)$
exists and satisfies $|f'(w)|\leq 1$.}
\vspace{.1in}

Angular limit means that $z$ is restricted to any angle
$\epsilon<\arg(z-w)<\pi-\epsilon$ if $w\in\bbR $,
or $\epsilon<\arg z<\pi-\epsilon$ if $w=\infty$.

Angular derivative is
the angular limit
$$f'(w)=\lim_{z\to w} (f(z)-f(w))/(z-w)$$
if $w\in\bbR $;
if $w=\infty$ then it is defined by the angular limit
$$\frac{1}{f'(\infty)}=\lim_{z\to\infty} f(z)/z.$$
The point $w$ in this theorem is called the Denjoy--Wolff
point of $f$. If $w\in\R\cup\{\infty\}$ is such a point
that the angular limit $\lim_{z\to w}f(z)=w$ and
the angular derivative $|f'(w)|\leq 1$, then $w$ is
the Denjoy--Wolff point.

Proof of Theorem~6. 

For the equations
$F_n(z)=0$, $G_n(z)=0$, and $H_n(z)=0$,
we will show that each of them
has at most one solution in the upper half-plane.
This solution belongs to
the second quadrant.

The Laplace transform of $k(t)$ is
\begin{eqnarray*}
K(z)&=&\int_0^\infty e^{-zx}\int_0^\infty e^{-tx}d\mu(t)dx\\
    &=&\int_0^\infty\int_0^\infty e^{-x(z+t)}dx\, d\mu(t)
    =\int_0^\infty\frac{d\mu(t)}{z+t}.
\end{eqnarray*}
This is called the Cauchy transform of the measure $d\mu$.
Condition (\ref{1}) ensures that the integral defining
$K$ is absolutely and uniformly
convergent on every compact in the $z$-plane that does not
intersect the negative ray. So $K$ is analytic
in the plane minus the negative ray.

Moreover,
\begin{equation}\label{Im}
\Ima K(z)\Ima z< 0,\quad z\in\bbC \backslash\bbR_-.
\end{equation}

We rewrite the equation $F_n(z)=0$
as
$z=f(z):=-n^2K(z)$, then $f$ maps $\bbC_+$ to $\bbC_+$ and by Schwarz's
Lemma has at most one fixed point in the upper half-plane.
If $\Rea z\geq 0$, then
$$\Rea(z+n^2K(z))=\Rea z+n^2\int_0^\infty\frac{t+\Rea z}{|t+z|^2}
d\mu(t)>0,$$
so the solution must lie in the second quadrant.

For $G_n(z)=0$, we first prove that there
are no solutions in the first quadrant.
Indeed $z^2+an^2$ maps the first quadrant into $\bbC_+$,
while
$n^2K(z)$ has negative imaginary part in the first
quadrant.

To prove that $G_n(z)$ has at most one zero in
$\bbC_+$ we consider a branch $\phi$ of the square root
which maps the lower half-plane onto the second quadrant.
Then the equation is equivalent to
\EQ{Fn}{
z=f(z):=n\phi(K(z)-a),
}
because all solutions are in the second quadrant.
Function $f$ maps $\bbC_+$ into itself, and thus by Schwarz's
Lemma can have at most one fixed point in $\bbC_+$.

Similar argument applies to $H_n(z)=0$.
There
are no solution in the first quadrant.
Indeed, if
$$
z^2+\epsilon zn^2+n^2=n^2K(z)
$$
and $z$ is in the first quadrant, then the LHS is in $\bbC_+$ but the RHS is in
$\bbC_-$.
So the equation is equivalent to
$$
z=f(z):=n\phi(K(z)-\epsilon z-1),
$$
because all solutions are in the second quadrant.
Function $F$ maps $\bbC_+$ into itself, and thus by Schwarz's
Lemma can have at most one fixed point in $\bbC_+$.

This completes the proof of part (i) of
Theorem \ref{complex}.
\vspace{.1in}

To prove part (ii), we first notice that
\begin{equation}\label{prop1}
K(z)\to 0\quad\mbox{as}\quad z=re^{i\theta}, r\to\infty,
\end{equation}
uniformly with respect
to $\theta$ for $|\theta|<\pi-\dl$,
for any given $\dl>0$.
This will be expressed by saying that $K\to 0$ as
$z\to\infty$ {\em non-tangentially}.
To show this we use the following lemma.
\BE{lemma}
{\LB{order}
If
\EQ{dl1}{
|\arg z|<\pi-\dl,
}
then
\EQ{sum}{
|z+t|\asymp |z|+t,  \ t\ge0.
}
}

Proof of the lemma. First, $|z+t|\leq|z|+t.$
Second,
$$|z+t|=|z||1+t/z|\geq |z|\cos\delta,$$
and similarly
$$|z+t|=t|1+z/t|\geq t\cos\delta.$$
Thus $|z+t|\geq (1/2)(|z|+t)\cos\delta.$
This gives \eqref{sum}.
\vspace{.1in}

Now
$$
|K(z)|\le C\int_0^\iy \frac{d\mu(t)}{|z|+t},
$$
in the sector (\ref{dl1}),
and
we obtain \eqref{prop1}.
\vspace{.1in}

We rewrite the equation $G_n(z)=0$ as
$$
z_n=in\sqrt a\sqrt{1-K(z_n)/a}.
$$
As $K(z)\to 0$ by Lemma 9, we obtain (ii).
\vspace{.1in}

Now we prove (iii).
Condition \eqref2 permits to
obtain asymptotics of $K$:
\EQ{A}{
K(z)=A/z+o(|z|^{-1}),\quad z\to \infty,
}
uniformly with respect to $\arg z$
in $|\arg z|\leq \pi-\epsilon$. Indeed,
\begin{eqnarray*}
\l|K(z)-A/z\r| &=&\l| \int_0^\infty
\left(\frac{1}{z+t}-\frac{1}{z}\right)d\mu(t)\r|\\
&\le& \frac{1}{|z|}\int_0^\infty\frac{t}{|z|+t}d\mu(t)
=o(|z|^{-1}).
\end{eqnarray*}

Now we rewrite $F_n(z_n)=0$ as
$$
z_n=-n^2\left(A\frac{1}{z_n}+o(z^{-1})\right).
$$
which gives (iii).
This completes the proof of Theorem 6.
\vspace{.1in}

To find out whether a solution in $\bbC_+$ exists for a given $n$, we consider
the case of discrete measure $\mu$. In the case that
$\mu$ has a finite support,
$$K(z)=\sum_{k=1}^N\frac{a_k}{z+b_k},$$
our arguments are elementary.
In this case, $K(z)$ is a real rational function
with $N$ poles, all of them on the real line.
Then $F_n$ is a rational function
of degree $N+1$ because it has an
additional pole at infinity.
So it must have $N+1$ zeros in the complex plane.
On each interval $I_k=(-b_{k+1},-b_k)$, $1\leq k\leq N-1$
there is one zero by the Bolzano-Weierstrass Theorem.
The remaining two zeros can lie either one in
$\bbC_+$ and one in $\bbC_-$ or both on some interval $I_k$,
(which then contains three zeros),
or both on the interval $I_0=(-b_1,0).$
One can give examples of each possibility.
\vspace{.1in}

{\em Examples.

1. Suppose that the measure $d\mu$ has two atoms, that is
$$
k(t)=\frac1{10}\l(e^{-t}+e^{-2t}\r).
$$
Then it easy to check that
$$
F_1(z)=K(z)+z=\frac{10}{z+1}+\frac{1}{z+2}+z
$$
has 3 real zeros
(and no non-real zeros).
The additional two zeros are in $(-1,0)$.

2. Let
$$k(t)=e^{-t}+200e^{-50t}.$$
It is easy to check that the equation
$$
K(x)+x=\frac{1}{x+1}+\frac{200}{x+50}+x
$$
has $3$ roots on the interval $[-50,-1]$.}
\vspace{.1in}

In the case of a measure with finitely many atoms,
the question whether the Denjoy--Wolff point belongs to the
real line can be solved in finitely many steps
by using the criterion that all roots of a polynomial equation
are real, see, for example \cite{Gant}.

In the general case of a discrete measure $\mu$
with atoms at $b_k$, we denote $I_k=(-b_{k+1},-b_k),k\geq 1$,
and $I_0=(-b_1,0)$.
If there is a solution $w$ in the upper half-plane, it
must be the Denjoy--Wolff point of $f(z)=-n^2K(z)$.
If there is no solution in the upper half-plane,
then the Denjoy--Wolff point $w$ belongs to some interval
$I_k$, and that $-1\leq f'(w)\leq 1$.

So theoretically we can find out whether there
is a solution
in the upper half-plane, by iterating $f$,
starting from any
point in $\bbC_+$, for example from the point $z_0=i$.
The sequence
$z_k=f_n(z_{k-1})$, $k=0,1,\dots,$ must converge.
If it converges to a point in
$\bbC_+$, then this point is the unique solution in $\bbC_+$.
Convergence in this case is geometric.
If $z_k$ converges to a point on the real line
then there is no solution in $\bbC_+$, but convergence in
this case may be extremely slow, if $|f'(w)|=1$.
\vspace{.1in}

Let $w$ be the Denjoy--Wolff point of $f$.
If $w\in \bbC_+$, then $w$ is the unique solution
of $F_n(z)=0$
in $\bbC_+$. The case $w=\infty$
is excluded because the angular
limit of  $f(z)$ as $z\to\infty$ is $0$.
If $w\in\bbR $ and $w\in I_k$ for some $k$, then
$f'(w)\in[-1,1]$ and this $I_k$ is
the unique interval of the $I_j$
which contains two additional real
zeros of $F_n$.

We conclude that in the case $w\in \bbC_+$, each interval
$I_k,\; k\geq 1$ contains one solution of $F_n(z)=0$
while $I_0$ contains no solutions.
\vspace{.1in}

The situation with GP2 and KV are similar: it has a solution
in $\bbC_+$ if and only if the Denjoy--Wolff of the
function $f=n\phi(K-a)$ is in $\bbC_+$, and the Denjoy--Wolff
point can in principle be found by iteration.
\vspace{.1in}

Now we prove the theorem \ref{complex1}.

\BE{lemma}{\LB{Prop2}
Under the assumption (\ref{mu})
we have
$$
K(z)=
\frac{b\pi\rho}{\sin\pi\rho}z^{\rho-1}
+O(z^{\alpha-1}),\quad
|z|\to\infty,
$$
uniformly with respect to $\arg z$
in any angle $|\arg z|<\pi-\dl$.
Here the use the principal branch of $z^{\rho-1}$ which
is positive on the positive ray.}

\vspace{.1in}

{\em Proof.}
First,
$$\int_0^\infty\frac{dt^\rho}{z+t}=
\frac{\pi\rho}{\sin\pi\rho}z^{\rho-1},$$
see, for example, \cite[Probl. 28.22]{E}
or \cite[Probl. 878]{Volk}.

So it is sufficient to prove our Lemma
for the case that $\mu(t)=O(t^\alpha)$.
We integrate
by parts:
$$
K(z)=\lim_{R\to\iy}\int_0^R\frac{d\mu(t)}{t+z}=
\lim_{R\to\infty} \l(\frac{\mu(R)}{z+R}+\int_0^R\frac{\mu(t)}{(t+z)^2}\r).
$$
In a sector $|\arg z|<\pi-\delta$ this gives
$$
|K(z)|\leq C
\lim_{R\to\infty}
\l(\frac{\mu(R)}{|z|+R}+\int_0^R\frac{t^\a}{t^2+|z|^2}dt\r)
$$
$$
=C\int_0^\iy\frac{t^\a}{t^2+|z|^2}dt=C_1|z|^{\a-1}.
$$
This proves the lemma.
\vspace{.1in}

Now the proof of Parts (i) and (ii) follows
the scheme of the proofs of
Theorem \ref{complex},
Parts (i) and (ii).
\vspace{.1in}

Proof of Theorem \ref{complexKV}.

In view of Theorem~6,
it is enough to prove that all roots of $H_n(z)$ are
real for $n$ large enough.
We rewrite the equation $H_n(z)=0$ as
$$n^2f(z)=z^2,$$
where $$f(z)=K(z)-\epsilon z-1.$$ We know from Theorem~6 that
all solutions belong to the second quadrant. Function $f$ maps
the upper half-plane into the lower half-plane. Let $\phi$ be
the branch of the square root that maps the lower half-plane
onto the second quadrant. Now we rewrite our equation as
$n\phi(f(z))=z.$ Function $f$ is strictly decreasing on
$(-\infty,-d)$ and $f(x)\sim -\epsilon x$ as $x\to-\infty$. So
there is a point $r<-d$ such that $f(x)>0$ on $(-\infty,r]$. We
have $n\phi(f(x))\sim -n\sqrt{-\epsilon x}$ as $x\to-\infty$,
It follows that
\begin{equation}
\label{last1}
n\phi(f(x))>x\quad\mbox{for}\quad x<x_0,
\end{equation}
where  $x_0<0$ depends on $n,K,\epsilon$.

Now suppose that $n$ is so large that
\begin{equation}
\label{last2}
n\phi(f(r))<r.
\end{equation}
This will hold for $n$ large enough because
$\phi(f(r))<0$ as we established above.
Comparison of (\ref{last1}) with (\ref{last2}) shows that
there must be a point $r_0\in (-\infty,r)$
such that
$$n\phi(f(r_0))=r_0\quad\mbox{and}\quad
\frac{d}{dx}n\phi(f(x))
\vert_{x=r_0}\in (0,1].$$
This shows that $r_0$ is an attracting fixed
point of the function
$n\phi(f)$, and application of the Denjoy--Wolff theorem
completes the proof.
\vspace{.1in}

%
We finish with the following

\emph{Question.} Can one extend Theorem~8
to arbitrary measure $d\mu$
satisfying (\ref{1}) ?

\vspace{.2in}

{\em Purdue University

West Lafayette IN 47907-2067

USA

eremenko@math.purdue.edu
\vspace{.2in}

State Marine Technical University,

St. Petersburg, Russia

sergei.ivanov@pobox.spbu.ru}

\begin{thebibliography}{1}

\bibitem{CMD} C. M. Dafermos, {\em Asymptotic stability
in viscoelasticity},
Arch. Rational Mech. Anal., 37 (1970), 297-308.
\bibitem{E} M. A. Evgrafov, Collection of problems in the theory
of analytic functions, Moscow, ``Nauka'', 1972 (Russian).

\bibitem{Gant} F. R. Gantmakher, Theory of Matrices,
Moscow, Nauka, 1988. English translation by AMS 1998.

\bibitem{Sham} A. A. Gavrikov, S.A. Ivanov, D.Yu. Knyazkov,
V.A. Samarain , A.S. Shamaev, V. V. Vlasov,
{\em Spectral properties of composite media},
Contemporary Problems of Mathematics and Mechanics,
v.1, 2009, 142-159 (Russian).

\bibitem{GuPip} M. E. Gurtin, A. C. Pipkin
{\em A general theory of heat conduction with finite
wave speeds}.
Archive for Rational Mechanics and Analysis 1968; 32:113-126.

\bibitem{I}S. A. Ivanov,
``Wave type'' spectrum of the Gurtin-Pipkin
equation of the second order,
arXiv; arxiv.org/abs/1002.2831, 8 p.
\bibitem{IS}S. A. Ivanov, T. L. Sheronova,
Spectrum of the heat equation with memory,
arXiv;  arxiv.org/abs/0912.1818v1, 10p.

\bibitem{P05}
L. Pandolfi,{\em The controllability of the
Gurtin-Pipkin equation:
a cosine operator approach}.
Appl. Math. Optim. 52 (2005), no. 2, 143--165.
\bibitem{Shapiro} J. Shapiro,
Composition operators and classical
function theory, Springer, 1993.

\bibitem{IP}
Ivanov, S., Pandolfi, L.,
{\em Heat equation with memory:
lack of controllability to rest,}
Journal of Math. Analysis and Appl., pp. 1-11, 2009, Vol. 355.


\bibitem{EMV} F. M.  Vegni,
{\em Dissipativity of a condensed phase field
systems with memory},
Discrete and continuous dynamical systems,
Volume 9, Number 4, July 2003.
\bibitem{VW}

V. V. Vlasov, J. Wu,
{\em Solvability and Spectral Analysis of Abstract Hyperbolic Equations with
Delay}.
Funct. Differ. Equ. 16 (2009), no. 4, 751-768.
\bibitem{Volk} L. I. Volkovyskii, G. Lunz and I. Aramanovich,
Collection of problems in the theory of functions
of a complex variable, Moscow, Fizmatgiz, 1960;
there is an English translation published by Dover.
\bibitem{W} D. Widder, Laplace transform, Princeton UP, 1946.
\end{thebibliography}
\end{document}